\documentclass[12pt]{amsart}
\usepackage{amssymb,amscd,xy}
\xyoption{all}

\theoremstyle{plain}
\newtheorem*{thm}{Theorem}
\newtheorem*{lem}{Lemma}
\newtheorem*{prop}{Proposition}

\theoremstyle{remark}
\newtheorem*{rem}{Remark}
 
\newtheorem*{notat}{Notations and Conventions}

\def\Hom{\operatorname{Hom}}
\def\Ker{\operatorname{Ker}}

\def\Res{\operatorname{Res}}
\def\Spec{\operatorname{Spec}}
\def\Id{\operatorname{Id}}
\def\Im{\operatorname{Im}}
\def\Ind{\operatorname{Ind}}
\def\Pic{\operatorname{Pic}}
\def\rad{\operatorname{rad}}
\def\Max{\operatorname{Max}}
\def\tilK{\widetilde{K_0}}

\def\tr{\operatorname{tr}}
\def\rk{\operatorname{rank}}
\def\GL#1{\operatorname{GL}(#1)}
\def\E#1{\operatorname{E}(#1)}
\def\GLn#1{\operatorname{GL}_n(#1)}
\def\GLm#1{\operatorname{GL}_m(#1)}
\def\SL#1{\operatorname{SL}(#1)}
\def\Aut#1#2{\operatorname{Aut}_#1(#2)}
\def\U#1{\operatorname{U}(#1)}

\def\M{\mathfrak{M}}

\def\P{\mathbf{P}}
\def\Z{\mathbb{Z}}
\def\C{\mathcal{C}}
\def\R{\mathcal{R}}

\def\*{\!*\!}
\newcommand{\cat}[1]{\operatorname{\mathsf{#1}}}

\begin{document}


\title[$K_0$ of Invariant Rings]%
   {$K_0$ of Invariant Rings and Nonabelian $H^1$}
\author{Martin Lorenz}
\thanks{Research supported in part by NSF Grant DMS-9618521}
\keywords{ring of invariants, group action, nonabelian 1-cohomology, Grothendieck group,
 induction map, Galois action, Hilbert's Theorem 90}
\subjclass{13A50, 13D15, 16E20, 16W20, 18G50, 19A49}

\maketitle

\begin{center}
Department of Mathematics\\
Temple University\\
Philadelphia, PA 19122-6094\\
\texttt{lorenz@math.temple.edu}
\end{center}

\begin{abstract}
We give a description of the kernel of the induction map $K_0(R)\to
K_0(S)$, where $S$ is a commutative ring and $R=S^G$ is the ring of
invariants of the action of a finite group $G$ on $S$. The description
is in terms of $H^1(G,\GL S)$.
\end{abstract}

\bigskip


\section*{Introduction} \label{intro}

This article is concerned with the relationship between $K_0(S)$ and $K_0(R)$,
where $S$ is a commutative ring and $R=S^G$ denotes the subring of invariants 
under the action of a finite group $G$ on $S$. 

Specifically, working under the assumption that the trace 
$\tr:S\to R,\ s\mapsto\sum_{g\in G}s^g$, is surjective,
we shall study the kernel of the induction map
$$
\Ind^S_R=K_0(f): K_0(R)\to K_0(S)
$$
that is associated with the inclusion $f:R\hookrightarrow S$.  
We will describe
an embedding of $\Ker(\Ind^S_R)$ into the cohomology set $H^1(G,\GL S)$. 
Moreover, we will endow
$H^1(G,\GL S)$ with a natural commutative monoid structure,
essentially coming from the ``block diagonal" maps 
$\operatorname{GL}_n\times\operatorname{GL}_m\to\operatorname{GL}_{n+m}
\hookrightarrow\operatorname{GL}$, 
such that our embedding identifies
$\Ker(\Ind^S_R)$ with the group of units $\U{H^1(G,\GL S)}$.
To further describe this unit group,
we define $S_H$, for any subgroup $H$ of $G$, to be the 
factor of $S$ modulo
the intersection of all maximal ideals of $S$ whose inertia group 
contains $H$. Letting $\rho_H: H^1(G,\GL S)\to H^1(H,\GL {S_H})$ 
denote the map that is given by restriction from
$G$ to $H$ and the canonical map $S\twoheadrightarrow S_H$,
our main result reads as follows.

\begin{thm} $\Ker(\Ind^S_R)\cong \U{H^1(G,\GL S)} = \bigcap_H
\Ker\rho_H$, where $H$ ranges over
all cyclic (or, equivalently, all) subgroups of $G$. 
\end{thm}


The theorem follows via direct limits from a
corresponding result, with $\GLn{\,.\,}$ in place of $\GL{\,.\,}$,
for the kernel of the maps $\P_n(f):\P_n(R)\to\P_n(S)$,
where $\P_n(\,.\,)$ denotes the set of isomorphism classes
of f.g. projective modules of constant rank $n$.

As applications, we present a version of Hilbert's 
Theorem 90 for Galois actions on commutative rings and quickly derive the 
(known) structure
of the Picard groups of linear and multiplicative invariants
(\cite{K}, \cite{L}). Some open problems are also discussed.


\begin{notat}
Throughout this note,
\begin{tabbing}
  XXXXXX\=$R = S^G$\qquad\= \kill          
  \>$S$ \> will be a ring, assumed commutative from \S\ref{add}
	onwards,\\
  \>$G$ \> will be a finite group acting by automorphisms on $S$;\\
  \> \> the action will be written $s\mapsto s^g$, and \\
  \>$R = S^G$ \> will denote the ring of $G$-invariants in $S$.
\end{tabbing}
We make the standing hypothesis that the trace map
$\tr:S\to R,\ s\mapsto\sum_{g\in G}s^g$ is surjective or,
equivalently:
\begin{equation*}
\text{\sl There exists an element $x\in S$ with $\tr(x)=1$.}
\tag{*}
\end{equation*}
\end{notat}


\section{Nonabelian $H^1$} \label{h1}

\subsection{Definition} \label{basics}

We recall the definition of nonabelian $H^1$ following 
\cite{Se} \S I.5 (or \cite{BS}, \cite{Se2} p.~123ff), except that our group 
actions are on the right.

Let $X$ be a $G$-group, with $G$-action written as $x\mapsto x^g$. 
A \textit{(1-)cocycle}
is a map $d:G\to X$ satisfying
$$
d(gg')=d(g)^{g'}d(g')\qquad (g,g'\in G)\ .
$$
The set of 1-cocycles of $G$ in $X$ will be denoted by $Z^1(G,X)$. 
Two cocycles
$d,e\in Z^1(G,X)$ are called \textit{cohomologous} if there exists an
element $x\in X$ satisfying $d(g)=x^ge(g)x^{-1}$ for all $g\in G$. This 
defines
an equivalence relation on $Z^1(G,X)$. The set of equivalence classes is
$$
H^1(G,X)\ ,
$$
a pointed set with distinguished element the class of the \textit{unit cocycle}
$\mathbf{1}(g)=1$ for all $g\in G$. 

\begin{rem} $H^1(G,X)$ parametrizes the
conjugacy classes of complements of $X$ in the split extension $X\rtimes G$,
exactly as in the familiar special case where $X$ is abelian
(cf.~\cite{Ro}, 11.1.2, 11.1.3).
\end{rem}

\subsection{Examples} \label{examples}
(1) Suppose $G$ acts trivially on $X$. Then $Z^1(G,X)=\Hom(G,X)$ and
$$
H^1(G,X)=\Hom(G,X)/X\ ,
$$
with $X$ acting by conjugation on $\Hom(G,X)$. The distinguished class
consists of $\mathbf{1}$ alone. \medskip

(2) If $G=\langle g\rangle$ is cyclic of order $m$, then each cocycle
$d\in Z^1(G,X)$ is determined by the element $x=d(g)\in X$, and the eligible
elements of $X$ are precisely those satisfying the condition
$x^{g^{m-1}}x^{g^{m-2}}\dots x^gx=1$. Moreover, if $d,e\in Z^1(G,X)$
correspond in this manner to $x,y\in X$, respectively, then $d$ and $e$ 
are cohomologous precisely
if there exists $z\in X$ with $x=z^gyz^{-1}$. Thus, writing $\sim$ for
the equivalence relation on $X$ determined by this condition, we have
$$
H^1(G,X)\cong \{x\in X : x^{g^{m-1}}x^{g^{m-2}}\dots x^gx=1\}/\sim \ .
$$

(3) If the order of $X$ is finite and coprime to $|G|$ then $H^1(G,X)$
is trivial, by the uniqueness part of the Schur-Zassenhaus Theorem
(cf.~\cite{Ro}, 9.1.2).

(4) If $X$ is a linear algebraic group over an algebraically closed field
whose characteristic does not divide $|G|$ and $G$ acts algebraically on 
$X$ then $H^1(G,X)$ is finite. This is essentially due to A.~Weil
who explicitly dealt with the case of a trivial $G$-action on $X$ 
(\cite{We}, cf.~also \cite{Sl}). The general case is an easy consequence.


\subsection{Functoriality and direct limits} \label{dirlim}

Suppose that we are given a homomorphism of groups $\alpha: G\to G'$ and 
a $G'$-group $X'$. Then $X'$ can be
viewed as a $G$-group via $\alpha$. Any map of $G$-groups (i.e.,
any group homomorphism compatible with the $G$-actions) $f:X'\to X$ gives 
rise to a map $Z^1(G',X')\to Z^1(G,X), \ d\mapsto f\circ d\circ\alpha$, 
and this map passes down to a map of pointed sets
$$
(\alpha,f)^1_*: H^1(G',X')\to H^1(G,X) \ .
$$
In particular, we have the \textit{restriction maps} 
$H^1(G,X)\to H^1(H,X)$ for subgroups $H$ of $G$ and the 
\textit{inflation maps} $H^1(G/N,X^N)\to H^1(G,X)$ for normal 
subgroups $N$. Here, $X^N$ denotes the 
$N$-invariants in $X$. We will write $f^1_*$ for $(\Id_G,f)^1_*$.

The following easy lemma is surely well-known, but I am not aware of
a reference. (For commutative cohomology, see \cite{Br},
Prop.~(4.6) on p.~195.) 

\begin{lem}
Let $(X_n,f_{mn})$ be a direct system of $G$-groups and let
$X=\varinjlim X_n$ be the direct limit with its induced $G$-action
\textnormal{(cf.~\cite{Bou3}, p.~A I.117)}. If all $f_{mn}:X_n\to X_m\ 
(m\ge n)$ are
injective then $H^1(G,X)\cong\varinjlim H^1(G,X_n)$ (as pointed sets).
\end{lem}

\begin{proof} Put $\varphi_{mn}=(f_{mn})^1_*:H^1(G,X_n)\to H^1(G,X_m)$
and $\varphi_n=(f_n)^1_*$, where $f_n:X_n\to X$ is
the canonical map. Then the relations $f_m\circ f_{mn}=f_n$ entail
$\varphi_m\circ \varphi_{mn}=\varphi_n$, and so there is a unique map
$\varphi : \varinjlim H^1(G,X_n)\to H^1(G,X)$ with $\varphi\circ\psi_n=
\varphi_n$, where $\psi_n : H^1(G,X_n)\to\varinjlim H^1(G,X_n)$ is the
canonical map. We show that $\varphi$ is bijective; the fact that $\varphi$
respects distinguished elements is clear.

For surjectivity, let $d\in Z^1(G,X)$ be given. Since $G$ is finite, there
is an $n$ with $d(G)\subseteq f_n(X_n)$, and so the class of $d$ in $H^1(G,X)$
belongs to $\Im\varphi_n\subseteq\Im\varphi$.

As to injectivity, we must show that if $a,b\in H^1(G,X_n)$ satisfy
$\varphi_n(a)=\varphi_n(b)$ then there exists $m\ge n$ with $\varphi_{mn}(a)=
\varphi_{mn}(b)$ (cf.\ \cite{Bou2}, Prop.~6 on p.~E III.62). Say $a$ and $b$
are the classes of $d,e\in Z^1(G,X_n)$, respectively. Then $f_n\circ d$
and $f_n\circ e$ are cohomologous in $Z^1(G,X)$, i.e., there exists $x\in X$
with $(f_n\circ d)(g)=x^g(f_n\circ e)(g)x^{-1}$ for all $g\in G$. Now
$x\in f_m(X_m)$ for some $m\ge n$, say $x=f_m(y)$. Then
$$
(f_m\circ f_{mn}\circ d)(g)=f_m(y^g)(f_m\circ f_{mn}\circ e)(g)f_m(y^{-1})\ .
$$
Since all $f_{mn}$ are injective, so are the maps $f_n$ (\cite{Bou2},
Remarque 1 on p.~E III.63). Hence, $(f_{mn}\circ d)(g)=
y^g(f_{mn}\circ e)(g)y^{-1}$ holds for all $g\in G$ which shows that
$f_{mn}\circ d$ and $f_{mn}\circ e$ are cohomologous in $Z^1(G,X_m)$, as
required.
\end{proof}


\subsection{The monoid structure of $H^1(G,\GL S)$} \label{monoid}

Lemma \ref{dirlim} implies in particular that 
in  $H^1(G,\GL S)\cong\varinjlim H^1(G,\GLn S)$, where $G$ acts on  
$\GL S$ and on $\GLn S$ via its action on $S$.
Our goal here is to endow $H^1(G,\GL S)$ with the structure of a 
commutative monoid.

More generally, let $X$ be any $G$-stable subgroup of $\GL S$ containing
the matrices
$$
s_{m,n}\underset{\text{def}}{=} (-1)^{(m+1)n} \begin{pmatrix}
0_{n\times m} & -1_{n\times n} & & & \\
1_{m\times m} & 0_{m\times n}  & & & \\
& & 1 & & \\
& & & \ddots &
\end{pmatrix} \ .
$$
Note that $\det s_{m,n}=1$. (In fact, it is not hard to show that the 
subgroup of $\GL S$ that is generated by the matrices $s_{m,n}$ consists 
of all monomial matrices with entries $\pm 1$ and with
determinant 1.) Hence $s_{m,n}\in\E S$, the group generated by
the elementary matrices (\cite{Ba}, Prop.~1.6 on p.~226). Furthermore,
$$
s_{m,n}^{-1}
\begin{pmatrix}
        a_{n\times n} & & & & \\
        & b_{m\times m} & & & \\
        & & 1 & & \\
        & & & \ddots &
\end{pmatrix}
s_{m,n}\ \ = \ 
\begin{pmatrix}
        b_{m\times m} & & & & \\
        & a_{n\times n} & & & \\
        & & 1 & & \\
        & & & \ddots &
\end{pmatrix} \ .
$$
Thus, if
$x=\left(
\begin{smallmatrix}
x_{n\times n} & & & \\
& 1 & & \\
& & \ddots &
\end{smallmatrix}
\right)\in X_n=X\cap\GLn S$
then, for any $m$,
$$
x[m]\underset{\text{def}}{=} 
\begin{pmatrix}
        1_{m\times m} & & & & \\
        & x_{n\times n} & & & \\
        & & 1 & & \\
        & & & \ddots &
\end{pmatrix} 
= s_{m,n}^{-1}xs_{m,n}\in X_{m+n}\ .
$$

For $a,b\in H^1(G,X)$, define $a+b\in H^1(G,X)$ as follows. Choose
cocycles $d,e\in Z^1(G,X)$ representing $a$ and $b$, respectively. Since
$G$ is finite, we have $d\in Z^1(G,X_m)$ and $e\in Z^1(G,X_n)$ for suitable
$m$ and $n$. For $g\in G$, put
$$
(d\oplus_m e)(g) \underset{\text{def}}{=} d(g)\cdot e(g)[m] =
\begin{pmatrix}
        d(g)_{m\times m} & & & & \\
        & e(g)_{n\times n} & & & \\
        & & 1 & & \\
        & & & \ddots &
\end{pmatrix}\in X\ .
$$
It is trivial to verify that $d\oplus_m e$ is a cocycle, and we define
$a+b$ to be its class in $H^1(G,X)$. To show that this is well-defined,
we first note that the class of $d\oplus_m e$ is independent of $m$, as long
as $d(G)\subseteq X_m$. Indeed, if $t\ge 0$ then
$$
(d\oplus_{m+t} e)(g)=\left(s_{t,n}[m]\right)^{-1}(d\oplus_m e)(g)s_{t,n}[m]\ ,
$$
and $s_{t,n}[m]=s_{t,n}[m]^g\in X$. Thus, $d\oplus_{m+t} e$ and
$d\oplus_m e$ are cohomologous.
Now suppose
that $a$ and $b$ are also represented by the cocycles $d'\in Z^1(G,X)$
and $e'\in Z^1(G,X)$, respectively. So there are matrices $x,y\in X$
with
$$
d'(g)=x^gd(g)x^{-1}\quad\text{and}\quad e'(g)=y^ge(g)y^{-1}\ .
$$
Fix $r$ so that $x,y,d(G),e(G)$ are all contained in $X_r$, and hence so
are $d'(G)$ and $e'(G)$. By the foregoing, it suffices to show that
$d\oplus_r e$ and $d'\oplus_r e'$ are cohomologous. But, putting
$z=x\cdot y[r]\in X$, we have 
$(d'\oplus_r e')(g)=z^g(d\oplus_r e)(g)z^{-1}$. 
Thus, $a+b\in H^1(G,X)$ is indeed well-defined. Commutativity
and associativity of + follow along similar lines, thereby making 
$H^1(G,X)$
a commutative monoid with neutral element the class of the unit cocycle.
--- For a different description of the monoid structure of $H^1(G,X)$ for
$X=\GL S$, see Lemma \ref{stab}. 


\subsection{Monoid maps} \label{maps}

Let $S'$ be another ring that is acted on by a group $G'$
and suppose that we are given a group homomorphism $\alpha:G\to G'$
and a $G$-equivariant ring map $\phi:S'\to S$, where $G$ acts on $S'$ via 
$\alpha$. Then we obtain a map of $G$-groups $\GL{\phi}:\GL{S'}\to\GL S$.
Thus, if $X'\subseteq\GL{S'}$ is a $G'$-stable subgroup containing the
matrices $s_{m,n}\in\GL{S'}$ of \S\ref{monoid} and if $X\subseteq\GL S$
is $G$-stable containing the image of $X'$ 
then we have a map of $G$-groups $f:X'\to X$. The map induced 
on cocycles $Z^1(G',X')\to Z^1(G,X)$ (cf.~\S\ref{dirlim}) is easily seen 
to respect 
the operations $\oplus_m$ of \S\ref{monoid}, and hence the map 
$(\alpha,f)^1_*: H^1(G',X')\to H^1(G,X)$ of \S\ref{dirlim} is actually a 
monoid map.  This applies in particular to restriction and induction maps
(with $\phi=\Id_S$). Moreover, we have the following

\begin{lem}
Let $X\subseteq\GL S$ be a $G$-stable subgroup with $X\supseteq\E S$ 
and let $f:X\to X^{\textnormal{ab}}=X/[X,X]$ denote the canonical map. Then 
$f^1_*: H^1(G,X)\to H^1(G,X^{\textnormal{ab}})$ is a monoid map, where 
$H^1(G,X^{\textnormal{ab}})$ has its usual group structure.
\end{lem}

\begin{proof} Consider $d,e\in Z^1(G,X)$, say
$d,e\in Z^1(G,X_m)$. Then, by \cite{Ba}, Prop.~(1.7) on p.~226,
$$
(d\oplus_m e)(g) \equiv  d(g)\cdot e(g) =
\begin{pmatrix}
        d(g)_{m\times m}e(g)_{m\times m} & & & \\
        & 1 & & \\
        & & \ddots &
\end{pmatrix} \mod{\E S}\ .
$$
Inasmuch as $X^{\text{ab}}=X/\E S$ (cf.~\cite{Ba}, Thm.~(2.1) on p.~228),
our assertion follows.
\end{proof}


\subsection{Units} \label{units}

The kernel of $\Ind^S_R: K_0(R)\to K_0(S)$ will turn out to be isomorphic
to the group of units $\U{H^1(G,\GL S)}$ of the monoid $H^1(G,\GL S)$
(see \S\ref{pfthm}). Here, we make some preliminary observations on the
unit group
$$
\U{H^1(G,X)}\ ,
$$
where $X$ is any $G$-invariant subgroup of $\GL S$ containing the matrices
$s_{m,n}$, as in \S\ref{monoid}.

\begin{lem}
Let $N=\Ker_G(X)$ denote the kernel of the action of $G$ on $X$. Then 
$\U{H^1(G,X)}\cong\U{H^1(G/N,X)}$ via inflation. In particular, if $G$ acts 
trivially on $X$, then $\U{H^1(G,X)}$ contains only the unit class.
\end{lem}

\begin{proof} Use the inflation-restriction sequence 
$H^1(G/N,X)\to H^1(G,X)\to H^1(N,X)$ (\cite{Se}, \S I.5.8).
This sequence is exact, the first map (inflation) is injective, and
both maps are monoid maps. Thus it induces an exact sequence of groups
$$
1\to \U{H^1(G/N,X)}\longrightarrow \U{H^1(G,X)}\longrightarrow 
\U{H^1(N,X)}\ .
$$
Part (ii) therefore reduces to the claim that $\U{H^1(N,X)}$ is trivial.
To verify this, recall that $H^1(N,X)=\Hom(N,X)/X$,
with $X$ acting by conjugation on $\Hom(N,X)$, and the unit class
consists of the unit map $\mathbf{1}$ alone
(cf.~\S\ref{examples}).
Thus, letting $\langle\,.\,\rangle$ denote $X$-conjugacy classes, the
equation $\langle d\rangle+\langle e\rangle=\langle \mathbf{1}\rangle$ for
$d,e\in\Hom(N,X)$ is equivalent with $(d\oplus_m e)(g)=1$ for all
$g\in N$, where $m$ is chosen as above. But the latter condition says
that $d=e=\mathbf{1}$.
\end{proof}


\section{The Kernel of Induction} \label{kerind}

\subsection{The skew group ring} \label{skew}

We will let 
$$
T=S\* G
$$ 
denote the \textit{skew group ring} that is associated
with the given $G$-action on $S$.
Thus $T$ is an associative ring containing $S$ as a subring and
$G$ is a subgroup of $\U T$, the group of units of $T$. The elements of $G$
form a free basis of $T$ as right $S$-module.
Multiplication in $T$ is based on the rule $ga\cdot hb = gha^hb$ for
$a,b \in S, g,h \in G$. The ring $S$ is an $R$-$T$-bimodule with action
$$
r\cdot a\cdot gb = ra^gb\qquad (r\in R,\ a,b\in S,\ g\in G)\ .
$$
Hypothesis (*) entails that $txt=t$, where we have put $t=\sum_{g\in G}g\in T$.
So $e=tx$ is an idempotent element of $T$ with
$eT=tT=tS\cong S_T$.
In particular, $S_T$ is projective and the ideal $I=TeT$ of $T$ satisfies
$I^2=I$ and $S_T\cdot I=S_T$.


\subsection{Some module categories} \label{categories}

Let $\cat{proj}R$ denote the category of finitely generated projective
(right) $R$-modules, and similarly for $T$, and let $\cat{add}S_T$ denote
the full subcategory of $\cat{proj}T$ whose objects are the direct summands
of the modules $S_T^n$ for $n\ge 0$.  The following lemma is well-known
but we include the proof for the reader's convenience.

\begin{lem}
\begin{enumerate}
\item[(i)] The functors $E:\cat{proj}R\to\cat{add}S_T$,
        $Q\mapsto Q\otimes_RS_T$ and
        $F:\cat{add}S_T\to\cat{proj}R$,
        $P\mapsto P^G$ yield an equivalence of categories
        $\cat{proj}R\approx\cat{add}S_T$.
\item[(ii)] A module $P$ in $\cat{proj}T$ belongs to $\cat{add}S_T$
        precisely if $PI=P$.
\end{enumerate}
\end{lem}

\begin{proof} (i) For $Q$ in $\cat{proj}R$, let
$\varphi_Q: Q\to (F\circ E)(Q)
= \left(Q\otimes_RS_T\right)^G$ denote the $R$-linear map given by
$\varphi_Q(q)=q\otimes 1$. Then $\varphi_R: R\to \left(R\otimes_RS_T\right)^G
\cong \left(S_T\right)^G=R$ is an isomorphism, and hence so is $\varphi_{R^n}$
for every $n$ and $\varphi_Q$ for every $Q$. Thus $\varphi$ is a
natural equivalence of functors $\Id_{\cat{proj}R}\cong F\circ E$.
Similarly, defining $\psi_P: (E\circ F)(P)=P^G\otimes_RS_T\to P$ for $P$ in
$\cat{add}S_T$ by $\psi_P(p\otimes s)=ps$, we obtain a natural equivalence
of functors $E\circ F\cong\Id_{\cat{add}S_T}$.

(ii) All modules $P$ in $\cat{add}S_T$ satisfy $PI=P$, because
$S\cdot I=S$.
Conversely, if $P$ in $\cat{proj}{T}$ satisfies $P=PI=PeT$ then, for some $n$,
$(eT)^n\cong S_T^n$ maps onto $P$, and so $P$ is a direct summand of $S_T^n$.
\end{proof}


\subsection{Another description of $\cat{add}S_T$} \label{add}

\textit{From now on, $S$ is assumed commutative.} We let $\Max S$
denote the set of maximal ideals of $S$. For
each $\M\in\Max S$, we put $G^Z(\M)=\{g\in G:\M^g=\M\}$, the
\textit{decomposition group} of $\M$, and
$$
T(\M)=(S/\M)\* G^Z(\M)\ ,
$$
the skew group ring that is associated with the action of $G^Z(\M)$ on
$S/\M$. As in \S \ref{skew}, $S/\M$ is a right module over $T(\M)$; this
module structure can be viewed as coming from $S_T$ by
restriction to $S\!* G^Z(\M)$ and reduction mod $\M$.
The following description of
\textbf{add}($S_T$) is adapted from \cite{Kr}, Proposition 3.

\begin{lem}
A module $P$ in $\cat{proj}T$ belongs to $\cat{add}S_T$ if and only
if, for all $\M\in\Max S$, there is an isomorphism of $T(\M)$-modules
$P/P\M\cong (S/\M)^r\ (r=\rk P_{\M})$.
\end{lem}

\begin{proof} The condition is surely necessary. For, if $P$ is a direct summand
of $S_T^n$, then $P/P\M$ is a direct summand of the homogeneous
$T(\M)$-module $(S/\M)^n$.

For the converse, consider some $\M\in\Max S$ and put
$\M^0=\bigcap_{g\in G}\M^g$, a $G$-stable ideal of $S$. Then
$$
S/\M^0\cong\bigoplus_{g\in G^Z(\M)\backslash G}S/\M^g\cong
(S/\M)\otimes_{T(\M)}T
$$
as $T$-modules. Similarly, $P/P\M^0\cong P\otimes_SS/\M^0
\cong (P/P\M)\otimes_{T(\M)}T$ as $T$-modules, with $G$ acting
``diagonally" on
$P\otimes_SS/\M^0$: $(p\otimes\bar{s})g=pg\otimes\bar{s}^g$.
 By hypothesis, $P/P\M$ is isomorphic to
$(S/\M)^r$ as $T(\M)$-modules, and so
$$
P/P\M^0\cong (S/\M^0)^r
$$
as $T$-modules. Since $S=S\cdot I$ (cf.~\S\ref{skew}), this isomorphism 
implies $P=PI+P\M^0$,
and since $\M$ was arbitrary, we further conclude that $P=PI$ 
(cf.~\cite{Bou}, Prop.~11 on p.~113). In view of Lemma \ref{categories}(ii),
this shows that 
$P$ belongs to $\cat{add}S_T$.
\end{proof}


\subsection{The induction triangle} \label{triangle}

For each $n\ge 0$, we let $\P_n(R)$ denote the set of isomorphism classes
of f.g. projective $R$-modules of constant rank $n$, and similarly for
$\P_n(S)$. These are pointed sets with distinguished elements
$\langle R^n\rangle$ and $\langle S^n\rangle$, respectively, where
$\langle\,.\,\rangle$ denotes isomorphism classes. Furthermore, $\P_{S,n}(T)$
will denote the set of isomorphism classes of f.g. projective right $T$-modules
having constant rank $n$ as $S$-modules, with distinguished
element $\langle S_T^n\rangle$. We have a commutative diagram of pointed sets
(cf.\ \cite{Ba}, Prop.~(7.3) on p.~130)
$$
\xymatrix{
\P_n(R)\
\ar@{^{(}->}[dr]_{\Phi_n=(\,.\,)\otimes_RS_T} \ar[rr]^{\P_n(f)}
    & & \P_n(S)\\
& \P_{S,n}(T) \ar[ur]_{\Res^T_{S,n}}
}
$$
By Lemma \ref{categories}(i), $\Phi_n$ is injective. The kernels
of the other two maps will be described in \S\S \ref{res} and \ref{pnf}
below. Recall that the \textit{kernel} of a map of pointed sets is defined
to be the preimage of the distinguished element of the target set.


\subsection{The kernel of $\Res^T_{S,n}$} \label{res}

We now consider the restriction map $\Res^T_{S,n}: \P_{S,n}(T)\to \P_n(S)$,
as in \S \ref{triangle}.  

\begin{lem}        
$\Ker(\Res^T_{S,n})\cong H^1(G,\GLn S)$ as pointed sets.
\end{lem}

\begin{proof} Each cocycle $d\in Z^1(G,\GLn S)$ gives rise to a $T$-module
structure $\left(S^n\right)_d$ on $S^n$ via
$$
x\cdot gs = x^gd(g)s\qquad (x\in S^n, g\in G, s\in S)\ .
$$
This action extends the regular $S$-module structure on $S^n$.
Conversely, if $\cdot$ is any right $T$-module operation on $S^n$ extending
the regular $S$-module structure then write, for $g\in G$,
$$
e_i\cdot g = \sum_{j=1}^ne_jd_{i,j}(g)\ ,
$$
where $e_i\in S^n$ is the basis element with 1 in the $i$-th position and
0s elsewhere and $d_{i,j}(g)\in S$. A routine verification shows that
$d=\left(d_{i,j}\right)$ is a cocycle of $G$ in $\GLn S$ and the given
$T$-module structure on $S^n$ is identical with $\left(S^n\right)_d$.

Since $\left(S^n\right)_d$  for the unit cocycle $d=\mathbf{1}$ is just 
$S_T^n$, we obtain a surjective map of pointed sets
$$
Z^1(G,\GLn S)\twoheadrightarrow \Ker(\Res^T_{S,n}),\quad d\mapsto
\langle\left(S^n\right)_d\rangle\ .
$$
Finally, for $d,e\in Z^1(G,\GLn S)$, we have
$\left(S^n\right)_d\cong\left(S^n\right)_e$ as $T$-modules precisely if
there is an $S$-module isomorphism 
$\left(S^n\right)_d\stackrel{\cong}{\longrightarrow}
\left(S^n\right)_e$ that commutes with the $G$-actions, that is, for some 
matrix $a\in\GLn S$, $(x\cdot g)a=(xa)\cdot g$ holds for all $x\in S^n,
g\in G$. The latter condition is equivalent with $x^gd(g)=x^ga^ge(g)a^{-1}$
which in turn just says that $d$ and $e$ are cohomologous. This completes 
the proof of the lemma.
\end{proof}

\begin{rem}
With $T(\M)$ as in \S\ref{add}, we have a map of pointed sets 
$\P_{S,n}(T)\to\P_{S/\M,n}(T(\M))$, $P\mapsto P/P\M$, which restricts to a 
map $\Ker(\Res^T_{S,n})\to\Ker\left(\Res^{T(\M)}_{{S/\M},n}\right)$. 
In terms of the identification
provided by the above Lemma, the latter becomes the map
$$
\rho'_{\M,n}: H^1(G,\GLn S)\to H^1(G^Z(\M),\GLn {S/\M})
$$
that is given by restriction from $G$ to $G^Z(\M)$ and reduction 
modulo $\M$.
\end{rem}


\subsection{The kernel of $\P_n(f)$} \label{pnf}

For each subgroup $H\le G$, we put
$$
J(H)=\bigcap_{\substack{\M\in\Max S\\ H\subseteq G^T(\M)}} \M \qquad
\text{and}\qquad S_H=S/J(H)\ ,
$$
where $G^T(\M)=\{g\in G : s^g-s\in\M\text{ for all }s\in S\}$ is the
\textit{inertia group} of $\M$. In addition to the maps
$\rho'_{\M,n}$ introduced in 
Remark \ref{res}, we will consider the analogous restric\-tion-reduction 
maps
$$
\rho_{\M,n}:H^1(G,\GLn S)\to H^1(G^T(\M),\GLn {S/\M})
$$
and
$$
\rho_{H,n}:H^1(G,\GLn S)\to H^1(H,\GLn {S_H})\ .
$$
Since the actions of $G^T(\M)$ on $\GLn {S/\M}$ and of $H$ on $\GLn {S_H}$
are trivial, we have 
$H^1(G^T(\M),\GLn {S/\M})=\Hom(G^T(\M),\GLn {S/\M})/\GLn {S/\M}$ 
and similarly for $H^1(H,\GLn {S_H})$
(cf.~\S\ref{examples}).
We let $\C$ denote the set of cyclic subgroups of $G$. 

\begin{prop}
As pointed sets,
\begin{alignat*}{2}
\Ker\P_n(f) &\cong \bigcap_{\M\in\Max S} \Ker\rho'_{\M,n}
& & =\bigcap_{\M\in\Max S} \Ker\rho_{\M,n} \\
& =\bigcap_{C\in\C} \Ker\rho_{C,n}
& & =\bigcap_{H\le G} \Ker\rho_{H,n} \ .
\end{alignat*}
\end{prop}

\begin{proof} In view of the induction triangle in \S\ref{triangle},
$\Ker\P_n(f)\cong\Im\Phi_n\cap\Ker(\Res^T_{S,n})$. Furthermore, by virtue of
Lemmas \ref{categories} and \ref{add}, if $\langle P\rangle\in\P_{S,n}(T)$ 
then
$\langle P\rangle\in\Im\Phi_n$ iff $\langle P/P\M\rangle$ is the distinguished
element of $\P_{S/\M,n}(T(\M))$ for all $\M\in\Max S$. Therefore, by Remark
\ref{res},
$$
\Im\Phi_n\cap\Ker(\Res^T_{S,n}) \cong \bigcap_{\M\in\Max S} \Ker\rho'_{\M,n} \ ,
$$
which establishes the $\cong$ in the proposition.
Now $\rho_{\M,n}=\Res^{G^Z(\M)}_{G^T(\M)}\circ\rho'_{\M,n}$, where 
$$
\Res^{G^Z(\M)}_{G^T(\M)}: H^1(G^Z(\M),\GLn {S/\M})\to 
H^1(G^T(\M),\GLn {S/\M})
$$
is the restriction map. Since this map has trivial kernel,
by the generalized ``Theorem 90" (\cite{Se}, Lemme 1 on p.\ 129 and 
\S5.8(a)), we conclude that
$\Ker\rho'_{\M,n}=\Ker\rho_{\M,n}$ which proves the first equality.

Consider $a\in H^1(G,\GLn S)$, say $a$ is the class of
$d\in Z^1(G,\GLn S)$. Then, in view of \S\ref{examples}(1),
\begin{align*}
a\in\bigcap_{\M\in\Max S}\Ker\rho_{\M,n}
&\Longleftrightarrow
\forall\M\in\Max S\ \forall g\in G^T(\M): 
\ \ d(g)\equiv 1_{n\times n}\mod{\M} \\
&\Longleftrightarrow
\forall g\in G\ \forall \M \text{ with } \M\supseteq J(\langle g\rangle): 
\ \ d(g)\equiv 1_{n\times n}\mod{\M} \\
&\Longleftrightarrow
\forall g\in G: \ \ d(g)\equiv 1_{n\times n}\mod{J(\langle g\rangle)} \\
&\Longleftrightarrow
a\in\bigcap_{C\in\C}\Ker\rho_{C,n}\ ,
\end{align*}
proving the second equality. Finally, the proof of 
$\bigcap_\M \Ker\rho_{\M,n} =
\bigcap_H \Ker\rho_{H,n}$ is completely analogous, based on the
observation that $H\subseteq G^T(\M)$ if and only if $\M\supseteq J(H)$.
This completes the proof of the proposition.
\end{proof}


\subsection{Stabilization} \label{stab}

For $m\ge n$, we now consider the \textit{stabilization maps}
$\P_n(R)\to\P_m(R)$, $\langle P\rangle\mapsto\langle P\oplus R^{m-n}\rangle$,
the analogous map for $S$, and the map $\P_{S,n}(T)\to\P_{S,m}(T)$,
$\langle Q\rangle\mapsto\langle Q\oplus S_T^{m-n}\rangle$. These are maps
of pointed sets which are compatible with the maps $\Phi_n$, $\Res^T_{S,n}$
and $\P_n(f)$ in the induction triangle (\S\ref{triangle}). Thus we obtain
a commutative diagram
$$
\xymatrix{
\varinjlim\P_n(R)\
\ar@{^{(}->}[dr]_{\varphi=\varinjlim\Phi_n} \ar[rr]^{\varinjlim\P_n(f)}
    & & \varinjlim\P_n(S)\\
& \varinjlim\P_{S,n}(T) \ar[ur]_{r=\varinjlim\Res^T_{S,n}}
}
$$
with $\varphi$ injective, by \cite{Bou2}, Prop.~7 on p.~E III.64.
Explicitly (cf.\ \cite{Bou2}, p.~E III.61),
$$
\varinjlim\P_n(R) = \left.\biguplus_{n\ge 0} \P_n(R)\right/\sim \ ,
$$
where $\biguplus$ denotes the disjoint union and $x\sim y$ for
$x=\langle P\rangle\in\P_n(R)$,
$y=\langle Q\rangle\in\P_m(R)$ iff $P\oplus R^{t-n}\cong Q\oplus R^{t-m}$
for some $t\ge\max(m,n)$. In other words, 
$\langle P\rangle\sim\langle Q\rangle$ iff $P$ and $Q$ are stably 
isomorphic. Now, $\biguplus_{n\ge 0} \P_n(R)$
is a commutative monoid under $\oplus$
and the equivalence relation $\sim$ respects this structure. Thus,
$\varinjlim\P_n(R)$ becomes a commutative monoid with identity element
the stable isomorphism class of $\{0\}$, that is, the f.g. free modules.
Actually, $\varinjlim\P_n(R)$ is a group: If $\langle P\rangle\in\P_n(R)$ is
given then $P\oplus Q\cong R^r$ for suitable $Q$ and $r$, and hence
$\langle P\rangle\oplus\langle Q\rangle=\langle R^r\rangle
\sim\langle 0\rangle$. In fact, letting $\tilK(R)$ denote the kernel of
$$
\rk:K_0(R)\to H_0(R)= \{\text{continuous maps} \Spec R\to\Z\}\ ,
$$
as usual (cf.~\cite{Ba}, p.~459),
the map sending $\langle P\rangle\in\P_n(R)$ to $[P]-[R^n]\in\tilK(R)$ 
passes down to a homomorphism of groups
$$
\varinjlim\P_n(R)\to\tilK(R)
$$
which is easily seen to be an isomorphism (cf.~\cite{W}, Chap.~II, 
Lemma 2.3.1).

The foregoing is valid for \textit{any} commutative ring $R$, and so also
applies to $\varinjlim\P_n(S)$. Under the identification $\varinjlim\P_n
\cong\tilK$, the top map $\varinjlim\P_n(f)$ of the 
above triangle becomes
$$
\tilK(f)=\Ind^S_R\Bigg|_{\tilK(R)}:\tilK(R)\to\tilK(S)\ .
$$

Things are similar for $\varinjlim\P_{S,n}(T)$: Isomorphism classes
$\langle X\rangle\in\P_{S,n}(T)$ and $\langle Y\rangle\in\P_{S,m}(T)$
become identified precisely if $X\oplus S_T^{t-n}\cong Y\oplus S_T^{t-m}$
for some $t\ge\max(m,n)$. However, since the submonoid of
$\left(\biguplus_{n\ge 0} \P_{S,n}(T),\oplus\right)$ that is
generated by $\langle S_T\rangle$ need no longer be cofinal, 
$\varinjlim\P_{S,n}(T)$ is merely a commutative monoid. The 
maps $\varphi$ and $r$ of the above triangle are monoid maps and 
$\varphi$ is mono, as was pointed out earlier. 

Thus, summarizing, we have the following stabilized induction
triangle of commutative monoids and groups
$$
\xymatrix{
\tilK(R)\
\ar@{^{(}->}[dr]_{\varphi=(\,.\,)\otimes_RS_T} \ar[rr]^{\tilK(f)}
    & & \tilK(S)\\
& \varinjlim\P_{S,n}(T) \ar[ur]_{r}
}
$$

\begin{lem}
$\Ker r\cong H^1(G,\GL S)$ as commutative monoids.
\end{lem}

\begin{proof} Since direct limits commute with kernels
(\cite{Bou2}, Cor.~(ii) on p.~E III.65), 
$$
\Ker r = \varinjlim\left(\Ker\Res^T_{S,n}\right)\ .
$$
Next, we infer from Lemma \ref{res} (and its proof) that
the stabilization map $\Ker\Res^T_{S,n}\to\Ker\Res^T_{S,n}$ $(m\ge n)$ 
translates into
the map $H^1(G,\GLn S)\to H^1(G,\GLm S)$ that is induced by $\GLn S\to\GLm S$,
$a\mapsto \left(\begin{smallmatrix} a & 0 \\ 0 & \mathbf{1} \end{smallmatrix}
\right)$. Thus, Lemma \ref{dirlim} implies that
$\Ker r\cong H^1(G,\GL S)$, at least as pointed sets. The fact that this
$\cong$ does respect the additive structures is a consequence of the obvious
isomorphism (using the notations of \S\S\ref{monoid},\ref{res})
$$
\left(S^m\right)_d\oplus\left(S^n\right)_e
\cong\left(S^{m+n}\right)_{d\oplus_m e}
$$
for $d\in Z^1(G,\GLm S)$ and $e\in Z^1(G,\GLn S)$.
\end{proof}

\subsection{Proof of the main result} \label{pfthm}

For each subgroup $H\le G$, let
$$
\rho_H:H^1(G,\GL S)\to H^1(H,\GL {S_H})=\Hom(H,\GL {S_H})/\GL {S_H}
$$
be given by restriction from
$G$ to $H$ and reduction modulo $J(H)$; so $\rho_H=\varinjlim \rho_{H,n}$ 
(cf.~\S\ref{pnf}). 
Similarly, for $\M\in\Max S$, put
$$
\rho_{\M}=\varinjlim\rho_{\M,n}: H^1(G,\GL S)\to H^1(G^T(\M),\GL{S/\M})\ .
$$
By \ref{maps},
these maps are monoid maps, and hence their kernels are submonoids of 
$H^1(G,\GL S)$.
It is now a simple matter to prove the Theorem stated in the Introduction.
Recall that $\C$ denotes the set of cyclic subgroups of $G$.

\begin{thm} $\Ker K_0(f)\cong\U{H^1(G,\GL S)}$, 
an isomorphism of groups, and
$$
\U{H^1(G,\GL S)}= \bigcap_{C\in\C}\Ker\rho_C=
\bigcap_{H\le G}\Ker\rho_H=\bigcap_{\M\in\Max S}\Ker\rho_{\M}\ .
$$ 
If $S$ is Noetherian of Krull dimension $d$ and $n>d$ then
$\Ker K_0(f)\cong \Ker \P_n(f)$.
\end{thm}

\begin{proof} Recall that $K_0(\,.\,)$ decomposes naturally as
$K_0(\,.\,)=H_0(\,.\,)\oplus\tilK(\,.\,)$, with
$\tilK(\,.\,)\cong\varinjlim\P_n(\,.\,)$. If the ring in
question is Noetherian of Krull dimension $d$, then
$\varinjlim\P_n(\,.\,)\cong\P_m(\,.\,)$ holds for any $m>d$ (\cite{Ba2}).

Applying this to the inclusion $f:R\hookrightarrow S$, we obtain that
$K_0(f)=H_0(f)\oplus\tilK(f)$. Here, $H_0(f)$ is injective
(cf.\ \cite{Ba}, Lemma (3.1) on p.~459), and hence
$$
\Ker K_0(f) = \Ker \tilK(f)\cong\varinjlim\Ker\P_n(f)
$$
(again using the fact that direct limits commute with kernels).
Moreover, if $S$ is Noetherian of Krull dimension $d$ (and hence so is $R$, by
virtue of hypothesis (*) ), then $\Ker K_0(f)\cong\Ker\P_n(f)$ for $n>d$ which
establishes the last assertion of the theorem.

Writing the formula for $\Ker\P_n(f)$ in Proposition \ref{pnf} as
$$
\Ker\P_n(f)\cong\Ker\rho_n\ ,
$$
where $\rho_n=\{\rho_{C,n}\}:H^1(G,\GLn S)\to
\prod_{C\in\C} H^1(C,\GLn {S_C})$, the stabilization map
$\Ker\P_n(f)\to\Ker\P_m(f)\ (m\ge n)$ becomes the map on $H^1$ that is
induced by the maps $\operatorname{GL}_n\to\operatorname{GL}_m$,
$a\mapsto \left(\begin{smallmatrix} a & 0 \\ 0 & \mathbf{1} \end{smallmatrix}
\right)$, for $S$ and $S_C$. Thus, using Lemma \ref{dirlim} and the fact
that direct limits
commute with kernels and with finite direct products
(cf.\ \cite{Bou2}, \textit{loc.~cit.\/} and Prop.~10/Cor.\ on
p.~E III.67/8), we deduce that
$$
\Ker K_0(f)\cong\Ker\rho \ ,
$$
where $\rho=\{\rho_C\}:H^1(G,\GL S)\to\prod_{C\in\C} H^1(C,\GL {S_C})$.
The $\cong$ is additive, by Lemma \ref{stab}.
Since $\Ker K_0(f)$ is a group, its image in $H^1(G,\GL S)$ must
be contained in the unit group $\U{H^1(G,\GL S)}$.
On the other hand, we infer from
Lemma \ref{units} that $\U{H^1(G,\GL S)}$ is contained in
$\bigcap_{H\le G}\Ker\rho_H$, which in turn is clearly contained in
$\Ker\rho=\bigcap_{C\in\C}\Ker\rho_C$. Hence equality must hold throughout, 
that is,
$$
\bigcap_{C\in\C}\Ker\rho_C=
\bigcap_{H\le G}\Ker\rho_H=\U{H^1(G,\GL S)}\ . 
$$
Finally, the equality $\bigcap_{H\le G}\Ker\rho_H=\bigcap_{\M\in\Max S}\rho_{\M}$
is proved exactly as in the proof of Proposition \ref{pnf}, thereby
completing the proof of the theorem.
\end{proof}

%
%


\section{Applications and Problems} \label{applications}

\subsection{Galois actions} \label{galois}

The $G$-action on $S$ is \textit{Galois}, in the sense of Auslander and
Goldman \cite{AG}, if and only if $G^T(\M)$ is trivial for all
$\M\in\Max S$ (cf.\ \cite{CHR}, Theorem 1.3 on p.~4). In this case, 
hypothesis (*) is satisfied (\cite {CHR}, Lemma 1.6 om p.~7).
Thus, by Proposition \ref{pnf} and Theorem \ref{pfthm}, 
$$
\Ker\P_n(f)\cong H^1(G,\GLn S)\quad\text{and}\quad
\Ker K_0(f)\cong H^1(G,\GL S)\ .
$$
In particular, $H^1(G,\GL S)$ is a group with the operation of \ref{monoid}.
In fact:

\begin{prop} If the action of $G$ on $S$ is Galois then 
$\Ker K_0(f)\cong H^1(G,\GL S)$ is annihilated by a power of $|G|$.
For $S$ Noetherian of Krull dimension $d$, $|G|^d$ will do.
\end{prop}

\begin{proof} By \cite{CHR}, Lemma 4.1 on p.~13, $S$ is f.g.\ projective of 
constant rank equal to $|G|$ as $R$-module. Therefore,
$$
[S_R]-|G|[R]=[S_R]-|G|1_{K_0(R)}\in \tilK(R)
$$
and, moreover, the restriction map $\Res^S_R:K_0(S)\to K_0(R)$ is 
defined. The composite $\Res^S_R\circ K_0(f):K_0(R)\to K_0(R)$ is clearly 
multiplication with $[S_R]$. Hence,
$$
\Ker K_0(f)\subseteq\operatorname{ann}_{K_0(R)}([S_R])\ .
$$ 
Recall that $\tilK(R)$ is a nil ideal of
$K_0(R)$, and if $S$ (or, equivalently, $R$) is Noetherian of Krull 
dimension $d$, then $\tilK(R)^{d+1}=\{0\}$ (\cite{Ba}, pp.~477 and 473).
Furthermore, $\Ker K_0(f)\subseteq \tilK(R)$ (cf.~the proof of
Theorem \ref{pfthm}). Hence,
$\Ker K_0(f)\cdot([S_R]-|G|)^t=\{0\}$ for some $t$, with $t=d$ a possible
choice for $S$ Noetherian of Krull dimension $d$. Consequently, 
$\Ker K_0(f)\cdot |G|^t=\{0\}$.
\end{proof}

The proposition can be viewed as an extension of Hilbert's ``Theorem 90"
(cf.~\cite{Se}, Lemme 1 on p.\ 129)
to commutative rings. 

\subsection{ } \label{question}
I don't know if Proposition \ref{galois} generalizes to arbitrary 
$G$-actions:
\smallskip\newline

\noindent\textbf{Question.}  
Is $\Ker K_0(f)\cong \U{H^1(G,\GL S)}$ always  
annihilated by $|G|^d$ if $S$ is Noetherian of Krull dimension $d$?
\smallskip\newline
 
This is indeed so for $d=0$, in which case 
$K_0(R)=H_0(R)$ and $K_0(f)=H_0(f)$ is mono, and for $d=1$, where 
$K_0(R)=H_0(R)\oplus\Pic(R)$ and $\Ker K_0(f)=\Ker\Pic(f)$ is isomorphic
to a subgroup of $H^1(G,\U S)$ (cf.~\S\ref{pic} below).
By a routine direct limit argument, a positive answer to the above question
would imply that the kernel of the induction map $K_0(f):K_0(R)\to K_0(S)$
is always $|G|$-primary (i.e., every element is annihilated by a power
of $|G|$), for any commutative ring $S$.
Finally, I note that 
the dual statement for $G_0$ is known to hold (\cite{BrL1} or \cite{BrL2}): The
cokernel of the restriction map $G_0(S)\to G_0(R)$ is
annihilated by $|G|^{d+1}$. The proof
given in \cite{BrL1} results from an analysis of the so-called \textit{coniveau
filtration} of $G_0$. The key to the above problem
might very well be the \textit{Grothendieck $\gamma$-filtration} (\cite{SGA6}, 
\cite{FL})
$$
K_0(R)=F^0_\gamma K_0\supseteq\tilK(R)=F^1_\gamma K_0\supseteq\ldots\supseteq
F^{d+1}_\gamma K_0=0\ .
$$
The first two slices are $F^0_\gamma/F^1_\gamma=H_0(R)$ and 
$F^1_\gamma/F^2_\gamma =\Pic(R)$. 
Not much appears to be known about the higher slices.


\subsection{Picard groups} \label{pic}

For any commutative ring $\R$, the set $\P_1(\R)$ of isomorphism classes
of f.g. projective $\R$-modules of constant rank 1 forms a group under
$\otimes_{\R}$, with identity element the distinguished element
$\langle\R\rangle$. This group is the \textit{Picard group} of $\R$,
usually denoted $\Pic(\R)$ (cf.\ \cite{Ba}, p.~131ff).

Specializing Proposition \ref{pnf} to the case $n=1$ and letting
$\U {\,.\,} = \operatorname{GL}_1(\,.\,)$ denote groups of units, 
we obtain the following result (cf.\ \cite{Kr,DMV,L}):

\begin{prop} There is an isomorphism of groups
\begin{align*}
\Ker\Pic(f) & \cong
\bigcap_{C\in\C} \Ker\left(H^1(G,\U S)\to\Hom(C,\U {S_C})\right) \\
& =  \bigcap_{H\le G} \Ker\left(H^1(G,\U S)\to\Hom(H,\U {S_H})\right)
\end{align*}
and an exact sequence of commutative monoids
$$
1\to \Ker\sigma\longrightarrow\Ker K_0(f)\longrightarrow
\Ker\Pic(f)\to 1\ ,
$$
where $\sigma: H^1(G,\SL S)\to\prod_{C\in\C} H^1(C,\SL {S_C})$.
\end{prop}

\begin{proof}
The fact that the isomorphism of Proposition \ref{pnf} is an isomorphism 
of groups, not just of pointed sets, for $n=1$ is
a consequence of the identity of $T$-modules $S_d\otimes_SS_e\cong S_{de}$
for $d,e\in Z^1(G,\U S)$. The exact sequence is a consequence of
the exact sequence
$1\to \SL{\,.\,}\to\GL{\,.\,}\stackrel{\det}{\longrightarrow}\U {\,.\,}\to 1$
which is split by the canonical embedding $\U {\,.\,}=
\operatorname{GL}_1(\,.\,)\hookrightarrow \GL{\,.\,}$. 
Indeed, this sequence leads to a commutative diagram of pointed sets
(\cite{Se}, \S\S 5.4, 5.5)
$$
\xymatrix{
1 \ar[d] & 1 \ar[d] \\
H^1(G,\SL S) \ar[d] \ar[r]^(.42){\sigma} &
       \prod_{C\in\C} H^1(C,\SL {S_C}) \ar[d] \\
H^1(G,\GL S) \ar[d]^{\pi_1} \ar[r]^(.42){\rho} &
       \prod_{C\in\C} H^1(C,\GL {S_C}) \ar[d]_{\prod\pi_{1,C}} \\
H^1(G,\U S) \ar[d] \ar[r]^(.42){\rho_1} \ar@/^1pc/[u]^{\mu_1}
       &\prod_{C\in\C} H^1(C,\U {S_C}) \ar[d] \ar@/_1pc/[u]_{\prod\mu_{1,C}}\\
1 & 1
}
$$
with $\pi_1\circ\mu_1=\Id$ and $\pi_{1,C}\circ\mu_{1,C}=\Id$. Here,
$\rho$ and $\rho_1$ are the usual restriction-reduction maps, as in
the proof of Theorem \ref{pfthm}.
The diagram yields the exact sequence
$$
1\to \Ker\sigma\longrightarrow\Ker\rho\longrightarrow
\Ker\rho_1\to 1
$$
which is in fact a sequence of commutative monoids, by \S\ref{maps}.
Finally, $\Ker K_0(f)\cong\Ker\rho$ and $\Ker\Pic(f)\cong\Ker\rho_1$.
\end{proof}


\subsection{Linear actions} \label{linear}

Here, $S=S(V)$ is the symmetric algebra of a finite dimensional $k$-vector
space $V$ and $G$ is a subgroup of $\GL V=\Aut k V$. The $G$-action on $V$
extends uniquely to an action of $G$ on $S$ by $k$-algebra automorphisms.
Hypothesis (*) amounts to the requirement that $|G|^{-1}\in k$, which will
be assumed, and linear actions are never Galois. Both assertions follow
from the existence of an augmentation $\varepsilon:S\to k$ which is
$G$-invariant (i.e., $\varepsilon(s^g)=\varepsilon(s)$ holds for all
$s\in S$, $g\in G$); it is given by $\varepsilon(V)=\{0\}$.

\subsubsection{The factors $S_H$} \label{linearfactor}
We now describe the factors $S_H=S/J(H)$ of $S$ that were
introduced in \S\ref{pnf}. Fix a subgroup $H$ of $G$ and let $V(H)$
denote the subspace of $V$ that is generated by the elements $v-v^h$
$(v\in V, h\in H)$. Then $H\subseteq G^T(\M)$ iff $V(H)\subseteq\M$.
Thus, $J(H)$ is the intersection of all $\M\in\Max S$ with $V(H)\subseteq\M$.
Now $V(H)S$ is a prime ideal of $S$; in fact, $S/V(H)S\cong S(V_H)$, where
$V_H=V/V(H)$ is the vector space of $H$-coinvariants of $V$.
Moreover, since $kH$ is semisimple, we have $V=V^H\oplus V(H)$.
Therefore, $J(H)=V(H)S$ and
$$
S_H\cong S(V_H)\cong S(V^H)\ .
$$
The canonical map $S\twoheadrightarrow S_H$ is $H$-equivariantly
split by $S_H\cong S(V^H)\hookrightarrow S$.

\subsubsection{Picard group \textnormal{(\cite{K})}} \label{linearpic}
Here, $\Pic(S) = 1$ and so $\Ker\Pic(f)=\Pic(R)$. Also, $\U S=k^*$ and all 
$\U {S_H}=k^*$, and so the intersection in Proposition \ref{pic}
reduces to the intersection of the kernels of the restriction maps
$\Hom(G,k^*)\to\Hom(H,k^*)$ which is obviously trivial. Thus,
$$
\Pic(R) = 1\ .
$$

\subsubsection{A problem of Kraft} \label{kraft}
Recall that $\Pic(R)$ is an image of $\tilK(R)$ (cf.~\S\ref{question}).
It is an open question, raised by Kraft  (\cite{Kr}, Problem 5.1), 
whether in fact 
$\tilK(R)$ is trivial or, equivalently, $\Ker K_0(f)=\{0\}$. A positive
answer to Question \ref{question} would easily entail this
in characteristics $>0$, and for fixed-point-free
actions in characteristic 0. In the latter case, Kraft's problem has 
already been settled by 
a different method by Holland \cite{Ho}, at least for algebraically closed base
fields $k$. The following Proposition gives a cohomological reformulation of
Holland's result, and of a result of Gubeladze \cite{Gu}. 
We put $S_+=VS$ and, as usual, $\GL{S,S_+}$ denotes the kernel of
the reduction map $\GL S \to \GL{S/S_+}=\GL k$, and similarly for $\GLn{S,S_+}$.

\begin{prop}
Assume that $k$ is algebraically closed of characteristic 0. If $G$ acts
fixed-point-freely on $V$ then $H^1(G,\GL{S,S_+})$ is trivial.
If, in addition, $G$ is abelian (and hence actually cyclic) then
$H^1(G,\GLn{S,S_+})$ is trivial for all $n$.
\end{prop} 

\begin{proof}
By fixed-point-freeness, $V(H)=V$ holds for all non-identity subgroups
$H$ of $G$, and hence $J(H)=S_+$. Thus, in the notation of Theorem \ref{pfthm},
$\Ker\rho_G\subseteq\Ker\rho_H$ and so $\U{H^1(G,\GL S)}=\Ker\rho_G$.
Finally, by \cite{Se}, Prop.~38 on p.~49, the split exact sequence of 
$G$-groups $1\to\GL{S,S_+}\to\GL S\to\GL{S/S_+}=\GL k\to 1$ gives rise
to an exact sequence of pointed sets
$$
1\to H^1(G,\GL{S,S_+})\longrightarrow H^1(G,\GL S)
\stackrel{\rho_G}{\longrightarrow} H^1(G,\GL{S/S_+})\to 1\ .
$$
Inasmuch as $\Ker\rho_G\cong\Ker K_0(f)$ is trivial, by \cite{Ho},
triviality of $H^1(G,\GL{S,S_+})$ follows (and conversely).

If $G$ is abelian then $R$ is an affine normal semigroup algebra
(cf. \cite{H}), and hence all f.g.~projectives over $R$ are free,
by Gubeladze's Theorem \cite{Gu}. In other words, all maps $\P_n(f)$
have trivial kernel which in turn says that 
$H^1(G,\GLn{S,S_+})$ is trivial, exactly as above, using Proposition \ref{pnf}
instead of Theorem \ref{pfthm}.
\end{proof}


\subsection{Multiplicative actions} \label{mult}

In this case, $S=kA$ is the group algebra of a f.g. free abelian group $A$
and $G$ is a subgroup of $\GL A =\Aut {\Z} A$ acting on $S$ by
means of the unique extension of the natural $\GL A$-action on $A$.
Again, hypothesis (*) is equivalent to $|G|^{-1}\in k$ and multiplicative
actions are never Galois, because of the $G$-invariant augmentation
$\varepsilon:S\to k$ given by $\varepsilon(A)=\{1\}$.

\subsubsection{The factors $S_H$} \label{multfactor}
Fix a subgroup $H$ of $G$,
let $[A,H]$ denote the subgroup of $A$ that is generated by the
elements $a^{-1}a^h$ $(a\in A, h\in H)$, and let $\omega([A,H])S$ denote
the ideal of $S$ that is generated by the elements $a-1$ with $a\in [A,H]$.
So $\omega([A,H])S$ is the kernel of the canonical map of group algebras
$S=kA\twoheadrightarrow kA_H$, where $A_H=A/[A,H]$ denotes the
$H$-coinvariants of $A$. Clearly,
$H\subseteq G^T(\M)$ iff $\omega([A,H])S\subseteq\M$. We claim that
$\omega([A,H])S$ is a semiprime ideal of $S$. Since $|H|$ is nonzero in
$k$, this will follow if we can show that the torsion group of $A_H$ is
annihilated by a power of $|H|$. To this end, write
$(\,.\,)'=(\,.\,)\otimes_{\Z}\Z[1/|H|]$ and view $A'$ as a module over the
group ring $\Z'H$. Putting $e=|H|^{-1}\sum_{h\in H}h\in\Z'H$, an idempotent
of $\Z'H$, we have $A\subseteq A'=\left(A'\right)^e\oplus \left(A'\right)^{1-e}
=\left(A'\right)^H\oplus [A',H]$. Thus $A'/[A',H]$ is torsion-free, and
hence so is $A/(A\cap [A',H])$. Since every element of $[A',H]/[A,H]$ is
annihilated by a power of $|H|$, our claim follows. We conclude that
$J(H)=\omega([A,H])S$ and
$$
S_H\cong kA_H \ .
$$
The canonical map $S\twoheadrightarrow S_H$ is not split, in general, as
$S_H$ need not be a domain.

\subsubsection{Picard group \textnormal{(\cite{L})}}
\label{multpic}
Again, $\Pic(S) = 1$ and so $\Pic(R)$ can be determined from 
Proposition \ref{pic}.
Here, $\U S=k^*\times A$ and $\U {S_C}=k^*\times U_1$ with
$A_C\subseteq U_1$, the group of normalized (augmentation 1) units
of $S_C=kA_C$ (a strict inclusion, in general). The maps
$H^1(G,\U S)=\Hom(G,k^*)\times H^1(G,A)\to\Hom(C,\U {S_C})=
\Hom(C,k^*)\times\Hom(C,U_1)$ decompose as the direct product of the
restriction maps $\Hom(G,k^*)\to\Hom(C,k^*)$ times the maps
$H^1(G,A)\to H^1(C,A)\to H^1(C,A_C)=\Hom(C,A_C)\hookrightarrow
\Hom(C,U_1)$. The first factor contributes nothing to the kernel, as for 
linear actions. Since the map $H^1(C,A)\to H^1(C,A_C)$ is mono for cyclic
$C$ (e.g., \cite{B}, p.~79), the contribution from the second factor is
the kernel of the restriction map $H^1(G,A)\to H^1(C,A)$. Therefore,
$$
\Pic(R)=\bigcap_{C\in\C}\Ker\left(\Res^G_C:H^1(G,A)\to H^1(C,A)\right)\ .
$$
This group need not be trivial; the results of a systematic computer
aided investigation of all cases with $\rk A\le 4$ are reported in \cite{L}.


\subsection{Moding out the radical} \label{rad}

Returning to the general situation where $S$ is an arbitrary commutative 
ring satisfying (*), we briefly consider the reduction maps 
$p: \GL S\to \GL {S/\rad S}$ and $p_n: \GLn S\to \GLn {S/\rad S}$. Here,
$\rad S$ denotes the \textit{Jacobson radical} of $S$. The following lemma is
an application of Lemma \ref{res}. 

\begin{lem} The maps $(p_n)^1_*: H^1(G,\GLn S)\to H^1(G,\GLn {S/\rad S})$
have trivial kernel, and so does $p^1_*: H^1(G,\GL S)\to H^1(G,\GL {S/\rad S})$.
\end{lem}

\begin{proof} It suffices to consider the maps $(p_n)^1_*$; the case of $p^1_*$
then follows by taking $\varinjlim$. 

Using the identification $\Ker(\Res^T_{S,n})\cong H^1(G,\GLn S)$ of Lemma
\ref{res} and writing $\overline{S}=S/\rad S$ and $\overline{T}=\overline{S}\* G$
(cf.~\ref{skew}), the map $(p_n)^1_*$ becomes the map 
$\Ker(\Res^T_{S,n})\to\Ker(\Res^{\overline{T}}_{\overline{S},n})$ that is given
by $\langle P\rangle\mapsto\langle\overline{P}=P/P\rad S\rangle$. Say 
$\langle P\rangle$ belongs to the kernel of this map, that is, 
$\overline{P}\cong\overline{S}_{\overline{T}}^n$ or, equivalently,
$P/P\rad S\cong S_T^n/S_T^n\rad S$. Since $\rad S\subseteq\rad T$
(cf.~\cite{P}, Theorem 7.2.5), the Nakayama Lemma implies that $P\cong S_T^n$
(cf.~\cite{Ba}, Prop.~(2.12) on p.~90). Thus $\langle P\rangle$ is the
distinguished element of $\Ker(\Res^T_{S,n})$, as required.
\end{proof}


\end{document}